%% file: main.tex
\documentclass[11pt,a4paper]{article}
\usepackage{cames0}
\usepackage[utf8]{inputenc}

\usepackage{times}
\usepackage[linesnumbered,ruled,vlined]{algorithm2e} 
\usepackage{hyperref} 
\usepackage[numbers]{natbib}
\usepackage{graphicx}
\usepackage{float}
\usepackage{amsmath}
\usepackage{amssymb}
\usepackage{authblk}
\usepackage{setspace}

\title{\textbf{\Large Isogeometric Topology Optimization Based on Topological Derivatives}}
\author[1,*]{Guilherme Henrique TEIXEIRA}
\author[2]{Nepomuk KRENN}
\author[2]{Peter GANGL}
\author[1]{Benjamin MARUSSIG}

\affil[1]{\textit{Graz University of Technology, Institute of Applied Mechanics, Technikerstraße 4/II, 8010 Graz, Austria}}
\affil[2]{\textit{Johann Radon Institute for Computational and Applied Mathematics, Altenberger Straße 69, A-4040 Linz, Austria}}
\affil[*]{\small \textit{Corresponding Author e-mail: teixeira@tugraz.at}}

\date{} 
\begin{document}
\maketitle

\begin{abstract}
    Topology optimization is a valuable tool in engineering, facilitating the design of optimized structures. However, topological changes often require a remeshing step, which can become challenging.
    In this work, we propose an isogeometric approach to topology optimization driven by topological derivatives. The combination of a level-set method together with an immersed isogeometric framework allows seamless geometry updates without the necessity of remeshing. At the same time, topological derivatives provide topological modifications without the need to define initial holes \cite{Amstutz2006}.
    We investigate the influence of higher-degree basis functions in both the level-set representation and the approximation of the solution. 
    Two numerical examples demonstrate the proposed approach, showing that employing higher-degree basis functions for approximating the solution improves accuracy, while linear basis functions remain sufficient for the level-set function representation.
\end{abstract}

\textbf{Keywords:} topology optimization, isogeometric analysis,  topological derivative, level-set method, immersed methods, higher-degree basis function.

\input{section1_Intro}
%
\input{section2_ProblemDescription}

%
\input{section3_LevelSet}

%
\input{section4_TopologicalDerivatives}

%
\input{section5_CutElements}
%
\input{section6_NumericalResults}

%
\input{section7_Conclusion}

\section{Acknowledgement}
This work is supported by the joint DFG/FWF Collaborative Research Centre CREATOR (CRC -- TRR361/F90) at TU Darmstadt, TU Graz, RICAM and JKU Linz.

\bibliographystyle{plainnat}  
\bibliography{bibliography}     

\end{document}

%% file: section1_Intro.tex
\section{Introduction}
\begin{figure}[t]
    \centering
    \includegraphics[scale=1]{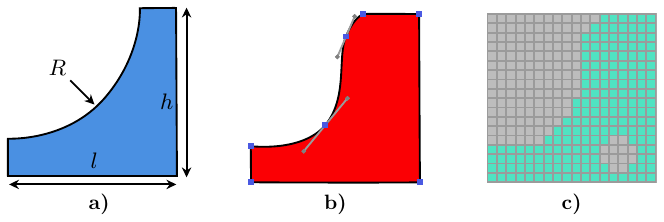}
    \caption{Different types of optimization: a) Parameter Optimization; b) Shape Optimization; c) Topology Optimization}
    \label{fig1}
\end{figure}

Design optimization describes an iterative process to define the optimal geometry of a structure given some constraints. This problem can be approached in several ways, including the optimization of geometric parameters, such as radius, length, or width (Figure \ref{fig1}a), the boundaries of the shape (Figure \ref{fig1}b) \cite{Wiesheu2024}, or the material distribution of the structure (Figure \ref{fig1}c) \cite{Ma2020}.
The last one is known as topology optimization, and since it was first introduced in \cite{Bendsøe1998}, several methods have been developed to approach the problem, and in consequence, various challenges have been addressed \cite{Sigmund2013}.
The most popular classes of topology optimization methods are density-based and level-set based methods.
One of the main challenges in level-set based topology optimization problems is related to a remeshing step, required due to the change in the domain during the optimization iterations.

The level-set method \cite{Osher1988}, extensively applied to shape and topology optimization \cite{Allaire2002, Allaire2004, Wang2003}, overcomes this situation by defining a level-set function to specify the region of a background mesh where material is defined while neglecting the contribution from the outside part. 
However, the level-set method still faces a problem with creating new holes. In the conventional level-set method \cite{Wang2003}, the evolution of the level-set is governed by the Hamilton-Jacobi equation. New holes cannot be directly introduced, and it depends on the holes of the initial geometry, which can then be merged or canceled in the optimization process. An alternative approach is to combine this equation with topological derivatives \cite{Burger2004, Allaire2006}, allowing the nucleation of new holes while the optimization is guided by the coupling of the Hamilton-Jacobi equation and the topological derivative, i.e., the sensitivity of the cost function with respect to pointwise material perturbations.
Another possible approach combining the level-set method and topological derivatives without solving the Hamilton-Jacobi is proposed in \cite{Amstutz2006} and then extended to multi-material in \cite{Gangl2020}. In this algorithm, the optimization is guided only by the topological derivative.

Isogeometric analysis (IGA), first introduced in \cite{Hughes2005}, presents the concept of connecting design and analysis using the same B-splines representing the geometry as basis functions. The straightforward control over the degree and smoothness of a B-spline basis is quite valuable for numerical simulations. Several research studies have been done using the isogeometric concept in variations of the conventional level-set method, which apply different discretizations for the level-set function, such as using radial basis functions \cite{Shojaee2012, Aminzadeh2022}, B-splines \cite{Jahangiry2017}, or piece-wise constant functions \cite{Khatibinia2020}. To the best of our knowledge, the combination of IGA and topological derivative-based level set optimization has only been considered in \cite{Roodsarabi2016}, where the conventional shape derivative-based level-set method extended by a topological derivative term was used.

In this work, we apply the approach of \cite{Amstutz2006} within the isogeometric framework, using B-splines for both the level-set function discretization and as basis functions to approximate the solution. The combination of the level-set method, topological derivatives, and IGA provides a simplified high-degree mesh, defined by the knot vector and control points. It also eliminates the need for remeshing and pre-definition of initial holes. In addition, the approach of \cite{Amstutz2006} does not require solving the Hamilton-Jacobi equation and depends only on the topological derivative. We study the sensitivity of the level-set representation in the optimized topology due to different polynomial degrees for approximating the level-set function and the solution. To accomplish this, we investigate two different settings: One with the same polynomial degree for both, and the other one with a linear level-set function discretization and a higher-degree approximation of the solution. 

Therefore, in Section 2, we present the linear elasticity problem investigated in the topology optimization approach and the main considerations for applying isogeometric analysis to it. Then, the definition of the level-set function and the derived topological derivative are provided, respectively, in Sections 3 and 4. In addition, Section 5 approaches the procedures to deal with the cut elements, and two numerical results applied to linear elasticity problems are shown in Section 6. Finally, in the conclusion, we summarize the main findings from the numerical examples.

%% file: section2_ProblemDescription.tex
\section{Problem Description}
 In this section, we present an overview of the problem to which our approach is applied. The problem studied is the compliance minimization in two-dimensional linear elasticity, and here we present the governing equation, the formulation of the minimization problem, and the key considerations for applying isogeometric analysis using an immersed approach based on the level-set method. 
\subsection{Linear Elasticity Problem} 
\begin{figure}[b]
    \centering
    \includegraphics[scale=1]{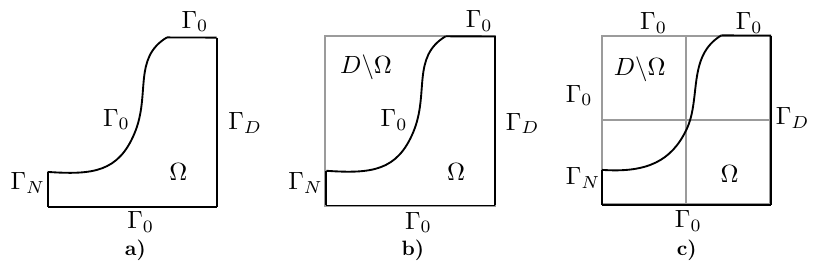}
    \caption{Representation of the domain problem: a) Domain $\Omega$ defined by the level-set; b) Domain $\Omega$ as a subset of domain $D$; c) Domain $\Omega$ inside of a B-spline background mesh defined from the knot vector of the geometry $D$ }
    \label{fig2}
\end{figure}
In this work, we consider linear elasticity problems defined on a domain $\Omega$, where the boundary $\partial\Omega$ is divided into three disjoint parts, such that $\partial\Omega=\Gamma_D\cup\Gamma_N\cup\Gamma_0$ and $\Gamma_D\cap\Gamma_N\cap\Gamma_0=\emptyset$, as shown in Figure \ref{fig2}a. In these three parts of the boundary, Dirichlet boundary conditions are applied in $\Gamma_D$, Neumann boundary conditions in $\Gamma_N$, and zero Neumann boundary conditions are applied in $\Gamma_0$. Therefore, the governing equations for the problem are given by
\begin{equation}\label{eq1}
    \begin{cases}
        -\nabla \cdot ( \boldsymbol{\sigma}(\mathbf{u})) = \mathbf{0} & \text{in } \Omega \\
        \mathbf{u} = \mathbf{0} & \text{on } \Gamma_D \\
        \boldsymbol{\sigma}(\mathbf{u}) \cdot \mathbf{n} = \boldsymbol{\tau} & \text{on } \Gamma_N \\
        \boldsymbol{\sigma}(\mathbf{u}) \cdot \mathbf{n} = \mathbf{0} & \text{on } \Gamma_0.
    \end{cases}
\end{equation}
Here, $\boldsymbol{\tau}$ is the load applied on the boundaries and $\mathbf{n}$ is the normal vector on it. While $\mathbf{u}$ represents the displacement field and $\boldsymbol{\sigma}$ is the stress tensor, which for linear elasticity and isotropic materials is defined as
\begin{equation*}
    \boldsymbol{\sigma(\mathbf{u})}=2\mu \boldsymbol{\epsilon}(\mathbf{u}) + \lambda\,\text{tr}(\boldsymbol{\epsilon}(\mathbf{u})) \mathbf{I},
\end{equation*}
where 
\begin{equation*}
\boldsymbol{\epsilon}(\mathbf{u})=\frac{1}{2} \left( \nabla \mathbf{u} + (\nabla \mathbf{u})^T \right)
\end{equation*}
is the strain tensor and, for 2d plain strain,
\begin{equation*}
 \mu = \frac{E}{2(1+\nu)}, \quad \lambda = \frac{E\nu}{(1+\nu)(1-2\nu)}
\end{equation*}
are the Lamé parameters, written with respect to the Young modulus $E$ and Poisson ratio $\nu$.

In this scenario, the domain $\Omega$ represents the material distribution of the geometry and is a subset of a larger domain $D$. This situation is graphically represented by Figure \ref{fig2}b.
\subsection{Immersed Isogeometric Approach}
The goal of the topology optimization is to find an optimal material distribution $\Omega$ under given constraints, such as boundary conditions or area penalization. This can be formulated as a minimization problem. In this scenario, the domain $\Omega$ changes during the optimization process, and solving the problem numerically would require redefining the mesh at each iteration.  To avoid the necessity of remeshing, we formulate the problem in the fixed domain D, as shown in Figure \ref{fig2}b, instead of $\Omega$, shown in Figure \ref{fig2}a, and we introduce a material property alpha, which is equal to $\alpha_{in}$ if it is inside of $\Omega$ and a small value $\alpha_{out}$ for outside. This approach is based on immersed methods, which are extensively applied to fluid mechanics, solid mechanics, interface problems, and several other areas. An extensive explanation of immersed methods and their aspects can be found in \cite{Wang2009, Prenter2023, Verzicco2023}. In this way, this approach allows us to define which part of $D$ represents $\Omega$, in such a way that the governing equation for the linear elasticity problem can be rewritten as
\begin{equation} \label{eq2}
     \begin{cases}
        -\nabla \cdot ( \alpha_\Omega\boldsymbol{\sigma}(\mathbf{u})) = \mathbf{0} & \text{in } D, \\
        \mathbf{u} = \mathbf{0} & \text{on } \Gamma_D, \\
        \boldsymbol{\sigma}(\mathbf{u}) \cdot \mathbf{n} = \boldsymbol{\tau}  & \text{on } \Gamma_N, \\
         \boldsymbol{\sigma}(\mathbf{u}) \cdot \mathbf{n} = \mathbf{0} & \text{on } \Gamma_0, 
     \end{cases}
    \quad\text{with}\quad
    \alpha_\Omega=
    \begin{cases}
        \alpha_{in} & \text{in } \Omega, \\
        \alpha_{out} & \text{on } D\setminus\overline{\Omega}, \\
    \end{cases}
\end{equation}
where $\alpha_{out}\ll1$ is a penalization parameter on the void, small enough to neglect the basis located outside the domain $\Omega$, but not too small to result in an ill-conditioned stiffness matrix \cite{Schillinger2015}.

Therefore, to finally solve the problem numerically, we discretize the domain $D$ using a background mesh, as shown in Figure \ref{fig2}c, where the basis functions used to approximate the solution field are B-splines of degree $p$, refined from the geometry $D$. These basis functions are constructed from a non-decreasing set of coordinates called knot vector
\begin{equation}
     \Xi=\left\{\xi_{1},\xi_{2},...,\xi_{n+p},\xi_{n+p+1}\right\}
\end{equation}
defined in a parameter space $\mathbb{P}=[\xi_1, \xi_{n+p+1}]$ of B-splines, where $n$ is the number of basis functions and $p$ is the polynomial degree. This construction is defined recursively starting from piece-wise constant functions for $p=0$
\begin{equation}
    \begin{aligned}
        B_{i,0}=\left\{
    \begin{matrix}
        1&\textnormal{if}\ \xi_i\le\xi<\xi_{i+1}\ \\
        0 & \text { otherwise } \\
        \end{matrix},\right.
    \end{aligned}
\end{equation}
and extended for higher degrees $p>0$ by applying the Cox-de Boor formula \cite{Hughes2005}
\begin{equation}
    B_{i,p}=\frac{\xi-\xi_i}{\xi_{i+p}-\xi_i}B_{i,p-1}(\xi)+\frac{\xi_{i+p+1}-\xi}{\xi_{i+p+1}-\xi_{i+1}}B_{i+1,p-1}(\xi).
\end{equation}
In sequence, having the basis functions in each parametric direction, the geometry mapping from the parametric domain $\mathbb{P}^2=[\xi_1, \xi_{n+p+1}]\times[\eta_1, \eta_{m+p+1}]$ to the physical domains $\mathbb{R}^2$ is then defined as
\begin{equation}
    \mathbf{x}(\xi,\eta) = \sum_{i=1}^{n}\sum_{j=1}^{m} B_{i,p}(\xi)B_{j,p}(\eta) \mathbf{C}_{i,j},
\end{equation}
where $\mathbf{C}_{i,j}$ are the control points that build the geometry.

Note that the basis functions used to approximate the solution field may have a different polynomial degree $p$ than those used to construct the geometry $D$. However, the geometric mapping remains based on the B-splines defined from the geometry $D$. Finally, the implementation of the problem is made within an open-source isogeometric analysis code \cite{Vázquez2016}, which provides the necessary features.
\subsection{Minimization Problem}
The minimization problem mentioned in the previous subsection, and which defines the topology optimization process, has the goal of searching for the optimal domain $\Omega\subset D$ that minimizes a cost function $J$. This expression can be written as
\begin{equation}
    \min_{\Omega \in \mathcal{E}} J(\Omega,\mathbf{u}), \\
\end{equation}
where $\mathcal{E}$ is a set of admissible subsets of $D$ and
\begin{equation}
    J(\Omega,\mathbf{u)}=\int_D \alpha_\Omega\boldsymbol{\sigma}(\mathbf{u}): \boldsymbol{\epsilon}(\mathbf{u}) \ dD + l\int_\Omega  \ d\Omega.
\end{equation}
Note that the area constraint, to avoid the solution to be $\Omega=D$, is addressed by the second term of the objective function and controlled by the parameter $l$ \cite{Krenn2021}.
%

%% file: section3_LevelSet.tex
\section{Discretized Level-set Representation}
\begin{figure}[t]
    \centering
    \includegraphics[scale=1]{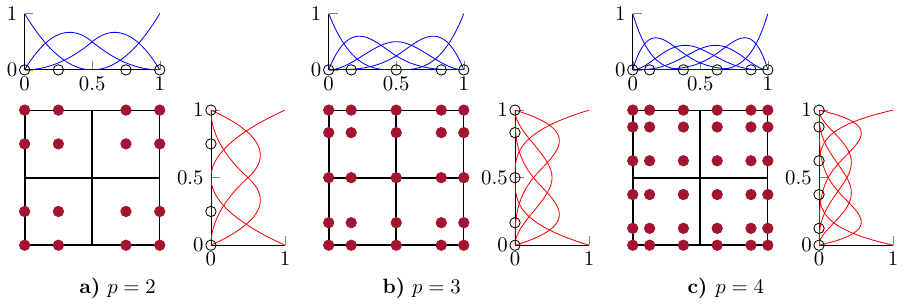}
    \caption{Distribution of the Greville abscissae on the elements for different polynomial degrees and basis functions defined by a) $\Xi=\{0 \:0 \: 0 \: 0.5 \: 1 \: 1 \: 1\}$. b) $\Xi=\{0 \: 0 \: 0 \: 0 \: 0.5 \: 1 \: 1 \: 1 \: 1\}$. c) $\Xi=\{0 \: 0 \: 0 \: 0 \: 0 \: 0.5 \: 1 \: 1 \: 1 \: 1 \: 1\}$ }
    \label{fig3}
\end{figure}
\begin{figure}[t]
    \centering
    \includegraphics[scale=1]{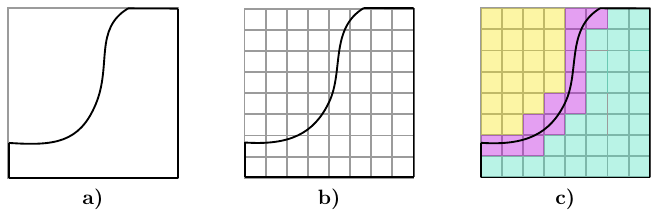}
    \caption{Type identification of the elements for assembling of the material property $\alpha$: a) Domain $D$ divided into two regions by a level-set function. b) Domains $D$ discretized as the background mesh. c) Identification of the elements. In yellow, the elements are located outside of $\Omega$. In blue, inside of $\Omega$ and in pink, the cut elements}
    \label{fig4}
\end{figure}
The domain $\Omega$ is represented by a continuous level-set function. This means that the interface that divides the material region $\Omega$ from the void region $D\setminus\Omega$ is defined by the zero set of the level-set function. Specifically, all points where the level-set function is smaller than 0 belong to $\Omega$, points where it is higher than zero belong to the void $D\setminus\Omega$, and points where the level-set function is equal to zero lie on the interface between the two regions
\begin{equation}
    \begin{cases}
        \phi(\mathbf{x}) < 0 & \iff \mathbf{x} \in \Omega, \\
        \phi(\mathbf{x}) = 0 & \iff \mathbf{x} \in \partial \Omega, \\
        \phi(\mathbf{x}) > 0 & \iff \mathbf{x} \in D \setminus \overline{\Omega}.
    \end{cases}
\end{equation}

 As the background mesh $D$ is discretized, we also discretize the domain $\Omega$. Therefore, the level-set discretization is made using B-spline basis functions of degree $d$, which might be the same or different from the degree $p$ of the basis functions used for approximating the solution,
\begin{equation}\label{eq10}
    \phi(\xi,\eta) = \sum_{i=1}^{n}\sum_{j=1}^{m} B_{i,d}(\xi)B_{j,d}(\eta) c_{i,j}.
\end{equation}

In sequence, the coefficients $c_{i,j}$ are obtained by solving a collocation problem, which enforces that the discretized level-set function (\ref{eq10}) is equal to the initial level-set function at the Greville abscissas. The position of the Greville abscissas, shown in Figure \ref{fig3}, works as anchors for the B-splines and is computed as
\begin{equation}
   \tilde{\xi}_{i}=\frac{\xi_{i+1}+\xi_{i+2}+...+\xi_{i+p}}{p}\ \ \ \ \ \ \ \ \ \ \ \ \ i=1,...,n.\ \ 
\end{equation}
During the optimization process, the evaluations of the level-set function at the Greville abscissas are updated, and the coefficients $c_{i,j}$ for the new level-set function discretization are obtained in the same process. 

In Figure \ref{fig4}, we observe the identification process of the region where the element is located. This is achieved by evaluating the level-set function at a group of points in each element. Then, based on the sign of the evaluations, we can identify the region of the element and attribute the corresponding material property or, in the cut element, compute the average of it based on the cut ratio of the element, which is given as follows
\begin{equation}
    \alpha|_T = \alpha_{out}+\frac{|T\cap\Omega|}{|T|}(\alpha_{in}-\alpha_{out}).
\end{equation}

%% file: section4_TopologicalDerivatives.tex
\section{Topological Derivatives}

In the previous sections, we defined the representation of a fixed domain $D$ and a level-set function $\phi$, used to represent the domain $\Omega$, both by B-spline discretization. In the following, we discuss how the topological derivative is computed and how the level-set function is updated, focusing on minimizing the cost function.

Considering $\mathbf{x}_0\in D\setminus\partial \Omega$, and defining $\omega_\varepsilon(\mathbf{x}_o)=\{\mathbf{x}\in\mathbb{R}^2:||\mathbf{x}_0-\mathbf{x}||<\varepsilon\}$ as a circular perturbation with radius $\varepsilon$ centered at $\mathbf{x}_0$. The introduction of the perturbation $\omega_\varepsilon(\mathbf{x}_o)$ in the domain $D$ results in a perturbed domain
\begin{equation}
    \Omega_\varepsilon =
    \begin{cases}
        \Omega \setminus \overline{\omega_\varepsilon(\mathbf{x}_0)} & \text{if } \mathbf{x}_0 \in \Omega, \\
        \Omega \cup \omega_\varepsilon(\mathbf{x}_0) & \text{if } \mathbf{x}_0 \in D \setminus \overline{\Omega}.
    \end{cases}
\end{equation}
Let $\mathcal J(\Omega) := J(\Omega, \mathbf u(\Omega))$ denote the reduced cost function where $\mathbf u(\Omega)$ denotes the unique solution to \eqref{eq1} for a given subdomain $\Omega$.
In this scenario, to measure the change in the cost function $J$ when a new hole around the point $x_0$ is introduced, the topological derivative is defined as
\begin{equation*} 
    d \mathcal J(\Omega)(\mathbf{x}_0)  := \lim_{\varepsilon \to 0} \frac{\mathcal J(\Omega_\varepsilon) - \mathcal J(\Omega)}{|\omega_\varepsilon|}= \lim_{\varepsilon \to 0} \frac{J(\Omega_\varepsilon,\mathbf{u}_\varepsilon) - J(\Omega,\mathbf{u})}{|\omega_\varepsilon|},
\end{equation*} 
where $\mathbf{u}_\varepsilon$ is the solution of (\ref{eq1}) replacing $\Omega$ by $\Omega_\varepsilon$.

To evaluate this expression, we adopt the approach proposed in \cite{GanglSturm2020}, which introduces the Lagrangian
\begin{equation*}
    L(\Omega,\mathbf{u},\boldsymbol{\lambda})=J(\Omega,\mathbf{u})+\boldsymbol{\boldsymbol{\lambda}} E(\Omega,\mathbf{u}),
\end{equation*}
where $E(\Omega,\mathbf{u})=0$ represents the weak form of the governing equation. This implies that  $L(\Omega,\mathbf{u},\boldsymbol{\boldsymbol{\lambda}})=J(\Omega,\mathbf{u})$ at the solution for all $\boldsymbol{\boldsymbol{\lambda}}$. 
Consequently, the topological derivative can be rewritten as
\begin{equation} \label{eq13}
    d \mathcal J(\Omega)(\mathbf{x}_0) = \lim_{\varepsilon \to 0} \frac{L(\Omega_\varepsilon,\mathbf{u}_\varepsilon,\boldsymbol{\boldsymbol{\lambda}}) - L(\Omega,\mathbf{u},\boldsymbol{\boldsymbol{\lambda}})}{|\omega_\varepsilon|}.
\end{equation} 

Plugging in the adjoint state $\boldsymbol{\boldsymbol{\lambda}}$  defined as the solution of $\partial_\mathbf{u}L(\Omega,\mathbf{u},\boldsymbol{\boldsymbol{\lambda}})=0$, and noting that $\boldsymbol{\boldsymbol{\lambda}}=-\frac12 \mathbf u$, after solving this limit for the linear elasticity problem, as shown in \cite{Amstutz2006sens}, an analytical expression is obtained, which depends only on the solution $\mathbf{u}$ and the material coefficient $\alpha$
\begin{equation}
    d \mathcal J(\Omega)(\mathbf{x}_0) =
\begin{cases}
     d \mathcal J_{in}(\Omega)(\mathbf{x}_0)=-3\alpha_{in}\left(\dfrac{\alpha_{out} - \alpha_{in}}{2\alpha_{out}+\alpha_{in}}\right)\boldsymbol{\sigma}(\mathbf{u}): \boldsymbol{\epsilon}(\mathbf{u}) - l & \text{if } \mathbf{x}_0 \in \Omega, \\[15pt] 
    d \mathcal J_{out}(\Omega)(\mathbf{x}_0)=-3\alpha_{out}\left(\dfrac{\alpha_{in} - \alpha_{out}}{2\alpha_{in}+\alpha_{out}}\right)\boldsymbol{\sigma}(\mathbf{u}): \boldsymbol{\epsilon}(\mathbf{u}) + l & \text{if } \mathbf{x}_0 \in D \setminus \overline{\Omega}.
\end{cases}
\end{equation}

From this, the generalized topological derivative is then defined as
\begin{equation}
    g_{\Omega}(\mathbf{x}) =
    \begin{cases}
    -d\mathcal J(\Omega)(\mathbf{x}) & \text{if } \mathbf{x} \in \Omega, \\
    d \mathcal J(\Omega)(\mathbf{x}) & \text{if } \mathbf{x} \in D \setminus \overline{\Omega},
    \end{cases}
\end{equation} 
and used to update the level-set, guiding the evolution of the domain $\Omega$. Algorithm 1 shows the update process. The update of the level-set is guided under a spherical linear interpolation, which uses the angle $\theta_i$, in $L^2$-sense, between the current level-set $\phi_{i}$ and the topological derivative $g_i$, as a parameter to define the next domain $\Omega_{i+1}$. Note that the stopping criterion is controlled by the same angle $\theta_i$, and this quantity works as a comparison between the current topological derivative $g_i$ and the level-set function $\phi_i$. Then, if $\theta_i=0$, the domain $\Omega_{i+1}$ is optimal and the topological derivative $g_i$ can be used as the level-set function $\phi$ \cite{Amstutz2006}. During this process, we apply a line search to define the parameter $\kappa$ used to update the level-set in the spherical linear interpolation. It is also applied some filtering processes, described in the section \ref{sec5}, to smooth the generalized topological derivative $g$, working similarly to a sensitivity filtering in density-based optimization \cite{Gangl2025}.
\begin{algorithm}[t]\label{alg1}
    \SetAlgoLined
    Initialize the level-set function $\phi_1 $\;
    \For{$i \gets 1$ \textbf{to} $n_{max}$}{
        Compute $g_{\Omega_i}(\mathbf{x}) =
        \begin{cases}
            -d \mathcal J(\Omega_i)(\mathbf{x}) & \text{if } \mathbf{x} \in \Omega_i, \\
             d \mathcal J(\Omega_i)(\mathbf{x}) & \text{if } \mathbf{x} \in D \setminus \overline{\Omega_i}.
        \end{cases}$ \;
         Compute $\theta_i=\arccos{(\frac{\langle \phi_i , g_{\Omega_i}\rangle}{||g_{\Omega_i}||_{L^2(D)}})}$\;
        \eIf{$\theta_i < \varepsilon_\theta$}{
            \textbf{break}\;
        }{
             $\phi_{i+1} = \frac{1}{\sin \theta_i} \left( \sin((1 - \kappa_i)\theta_i) \phi_i + \sin(\kappa_i \theta_i) g_{\Omega_i} \right)$\quad\\ where \quad $\kappa=\max\{1,\frac{1}{2},\frac{1}{4}, \dots\}$ such that $\mathcal J(\Omega_{i+1})<\mathcal J(\Omega_{i})$;
        }
        Update $c_{ij}$ in the discretization $\phi(\xi,\eta) = \sum_{i=1}^{n}\sum_{j=1}^{m} B_{i,d}(\xi)B_{j,d}(\eta) c_{ij} $
    }
    \Return $\phi$\;
            
    \caption{Level-set update}
\end{algorithm}
In the discretized setting, when the element $e$ of the background mesh $D$ is cut by the interface of the level-set function $\phi$, the generalized topological derivative is computed using a linear interpolation between the values computed inside and outside the domain $\Omega$, given by
\begin{equation}
    g_{\Omega}|_e(\mathbf{x}) =d \mathcal J_{out}(\Omega)(\mathbf{x})|_e+\frac{|e\cap\Omega|}{|e|}(-d \mathcal J_{in}(\Omega)(\mathbf{x})|_e-d \mathcal J_{out}(\Omega)(\mathbf{x})|_e).
\end{equation}

%% file: section5_CutElements.tex
\section{Cut Elements}\label{sec5}
    \begin{figure}[b]
        \centering
        \includegraphics[scale=1]{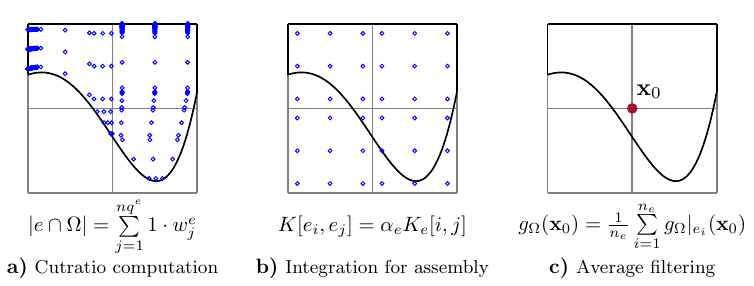}
        \caption{Approaches to treat the cut elements}
        \label{fig5}
    \end{figure}
The correct integration of the elements intersected by the level-set function and the precision on defining the proper material parameter in these elements play an important role in the quality of the results. Therefore, in this section, we present the procedures used to deal with the cut elements.
    
In the computation of the cut ratio $|e\cap\Omega| / |e|$, necessary to obtain the $\alpha$ property for the cut elements, we apply a quadrature library for implicitly defined geometries \cite{Saye2015, Saye2022}, which precisely follows the level-set function and provides quadrature points conforming with its zero isoline, as we can see in Figure \ref{fig5}a. This quadrature allows the integration of the regions defined by the zero level-set with high precision, also for complex geometries with high-degree representations of the interfaces. Therefore, with this quadrature, we can capture with precision the transition between the material property $\alpha$ in the inner part and the outer part of the domain $\Omega$. Some examples for area computation, moving geometries, and linear elasticity can be seen in \cite{Teixeira2025, toprak2025}.

However, when applying this precise quadrature rule for the assembly, we obtained some instabilities in the shape due to the integration of small parts of basis functions, located in the regions defined by the level-set function, and due to the jump between the material and the void regions, which results in a discontinuous field $\mathbf u$ and is approximated with higher-degree basis functions. To smooth these results, in the assembly process, we use the standard Gauss quadrature in the whole cut element and scale the local contribution by the corresponding material property $\alpha$, computed using the previous procedure, shown in Figure \ref{fig5}b. Therefore, this approach smooths out the discontinuity in the transition between the elements and results in a shape with less noise.

Another smoothing step may be applied for the generalized topological derivative. In particular, since we define an $\alpha$ property at each element, it occurs that,  when evaluating the derivative at a Greville point, shared by multiple elements, as we can see in Figure \ref{fig5}c, we have more than one $\alpha$ at the point. And to solve this, we consider the average of the derivative around that point. This procedure effectively creates a smooth transition where the material property changes and works well as a sensitivity filter \cite{Gangl2020}.

A second filtering to smooth the generalized topological derivative $g$ is performed by replacing $g$ in the spherical linear interpolation, as described in Algorithm 1, by the solution $ \tilde{g}$ of the PDE
\begin{equation}
    \begin{aligned}
        -\gamma&\Delta \tilde{g}_{\Omega} + \tilde{g}_{\Omega}=g_{\Omega} &\text{in } D\\
        &\nabla{\tilde{g}_{\Omega}}\cdot\mathbf{n}=0            &\text{on } \partial{D}
    \end{aligned}
\end{equation}

%% file: section6_NumericalResults.tex
\section{Numerical Results}
In this section, two numerical results for different geometries are presented. In both examples, the level-set is initialized as $\phi(\mathbf{x})=-1$, which results in a full material background mesh, and the background mesh is discretized with 128x128 elements. The material coefficients used are $\alpha_{in}=1$ for the material region $\Omega$ and $\alpha_{out}=10^{-4}$ for the void region $D\setminus \overline{\Omega}$. Additionally, the parameter for area control is $l=5$, the size control coefficient for the filtering is $\gamma=10^{-4}$, and the Young modulus and the Poisson ratio are 1 and 1/3, respectively. In these examples, we investigate the sensitivity of the level-set function representation by applying different polynomial degrees $d$ for the level-set function discretization and $p$ for approximating the solution. This is done considering two settings. One with the level-set function and the solution, both approximated with the same polynomial degree, and another with a linear level-set function representation and a higher-degree approximation of the solution. For the optimization algorithm, we consider a tolerance of $\varepsilon_\theta=1$ and a maximum number of iterations of 200.

\subsection{Cantilever}
    \begin{figure}[t]
        \centering
        \includegraphics[scale=1]{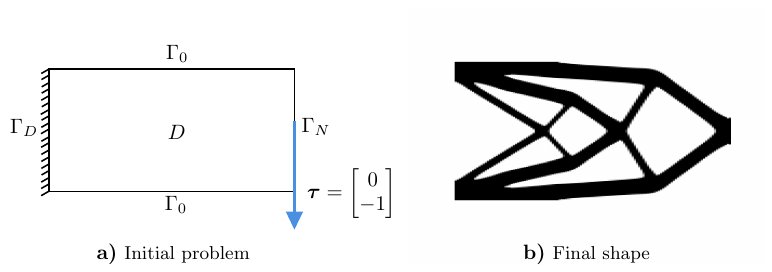}
        \caption{Cantilever problem}
        \label{fig6}
    \end{figure}
    \begin{figure}[t]
        \centering
        \includegraphics[scale=1.1]{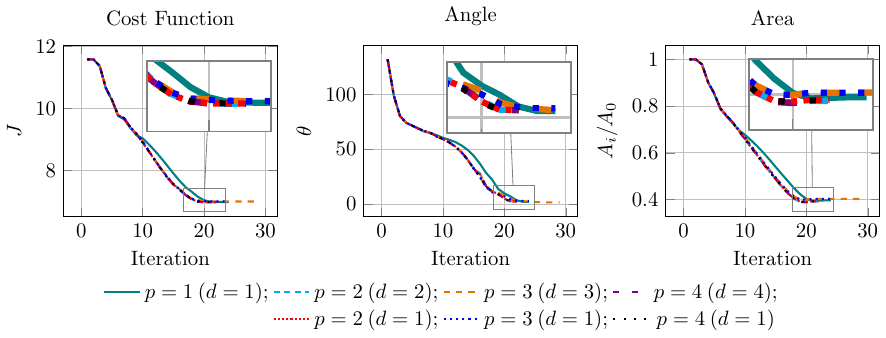}
        \caption{
        Comparison of the evolution of the cost function, angle, and area for different polynomial degrees of solution approximation, with a linear level-set function representation (dotted), and with a level-set function representation with the same polynomial degree as the solution approximation (dashed).}
        \label{fig7}
    \end{figure}
    \begin{figure}[t]
        \centering
        \includegraphics[scale=1.2]{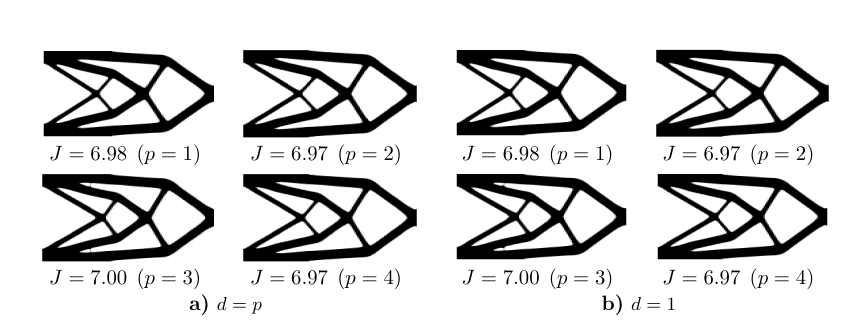}
        \caption{Final shape for different basis functions degree for approximating the solution: a) Level-set discretized with the same polynomial degree as the solution. b) Level-set discretized with linear basis functions}
        \label{fig8}
    \end{figure}
The first example is the cantilever problem, a benchmark example for topology optimization present in a large number of research papers, see, e.g., \cite{Allaire2002, Allaire2004, Amstutz2006, Shojaee2012, Roodsarabi2016, Jahangiry2017, Khatibinia2020, Krenn2021, Aminzadeh2022}. In our case, the domain $D$ is represented by the mapping from the parameter domain $\hat{D}=[0,1]\times[0,1]$ to a rectangle of size $2\times1$ with homogeneous Dirichlet boundary conditions on the left and a concentrated load on the right, as seen in Figure \ref{fig6}a. Figure \ref{fig6}b shows the final design for $p=d=2$, and Figure \ref{fig7} shows the evolution of the cost function, angle, and the area for the two settings. Additionally, the final shapes for different configurations of polynomial degree for the solution and the level-set function discretization are shown in Figure \ref{fig8}.

These results show that higher polynomial degrees $p$ provide a better convergence behavior compared to $p=1$, with a faster drop in the middle region of the graph. However, since we use a fine mesh ($128 \times 128$) to accurately represent the topology of the shape, increasing the polynomial degree $p$ does not necessarily increase the accuracy of the solution. In addition, in all cases, we converge to the same final solution, with similar convergence paths for higher degrees $p$, and the same number of steps is required for $p=2$ and $p=4$ in both level-set discretization settings. However, when using $p=3$, 29 steps are required for the higher-degree level-set discretization ($d = p$), compared to 24 for the linear level-set discretization ($d = 1$). We also observe a slight difference in the final part of the area graph between different solution discretizations using odd and even B-spline degrees. This difference happens because in the Greville abscissae for odd degrees, some points are shared between elements. As a result, an averaging is performed around these positions when evaluating the solution, effectively acting as a filtering process \cite{Gangl2020}. While this helps smooth the solution, it introduces a small difference compared to cases where such averaging is not required. In addition, looking at the final shapes in Figure \ref{fig8}, the linear representation of the level-set ($d=1$) produces a similar final shape compared to those obtained with the higher-degree discretizations ($d=p$).

\subsection{Quarter Ring}
\begin{figure}[t]
        \centering
        \includegraphics[scale=1]{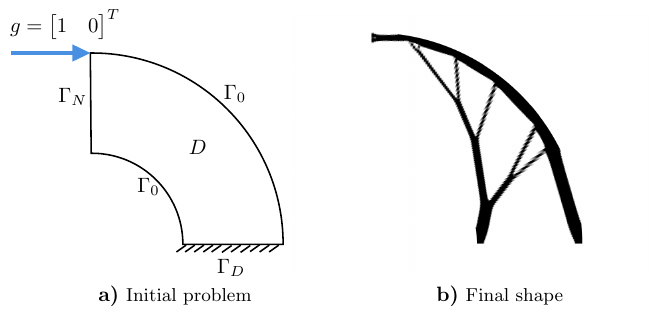}
        \caption{Curved cantilever problem}
        \label{fig9}
    \end{figure}
    \begin{figure}[t]
        \centering
        \includegraphics[scale=1.1]{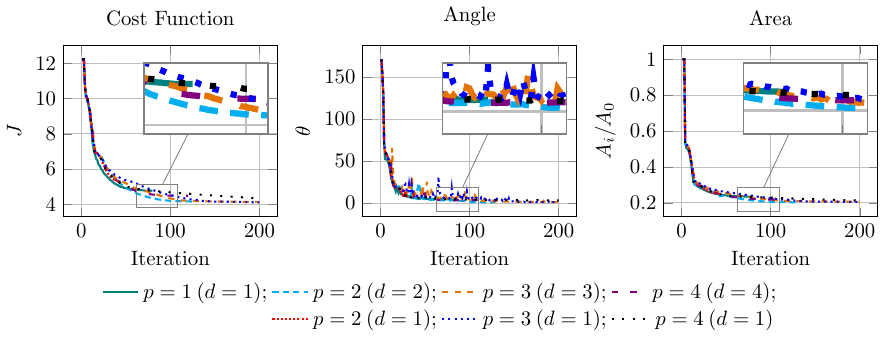}
        \caption{
        Comparison of the evolution of the cost function, angle, and area for different polynomial degrees of solution approximation, with a linear level-set function representation (dotted), and with a level-set function representation with the same polynomial degree as the solution approximation (dashed).}
        \label{fig10}
    \end{figure}
    \begin{figure}[t]
    \centering
    \includegraphics[scale=1.1]{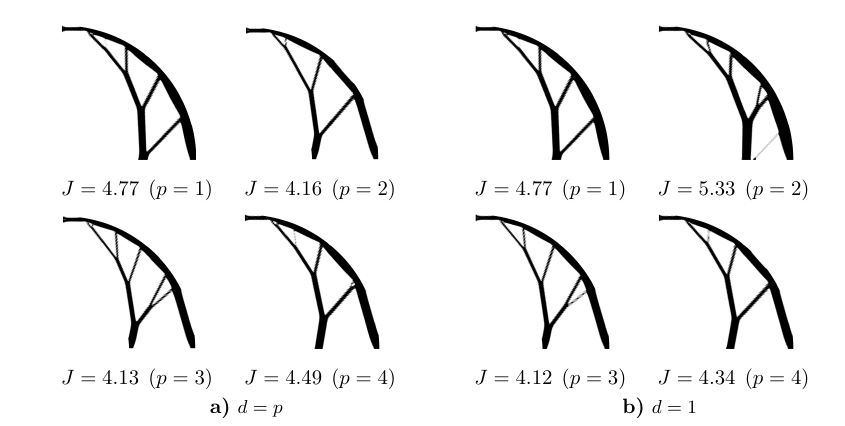}
    \caption{Final shape for different basis functions degree for approximating the solution: a) Level-set discretized with the same polynomial degree as the solution. b) Level-set discretized with linear basis functions.}
     \label{fig11}
\end{figure}
This example is also considered in several research papers \cite{dede2012, Jahangiry2017, Aminzadeh2022, Aminzadeh2024}, where the same geometry is applied under different approaches and loading configurations. In our example, the domain $D$ is defined by a mapping from the parameter domain $\hat{D}=[0,1]\times[0,1]$ to a quarter of a ring with inner radius $R_{in}=1$ and outer radius $R_{out}=2$. Homogeneous Dirichlet boundary conditions are imposed on the bottom boundary, and a concentrated load is applied at $(0,2)$, as illustrated in Figure \ref{fig9}a. The final design for $p=d=3$ is shown in Figure \ref{fig9}b, while the evolution of the level-set function and the corresponding final shapes are presented in Figure \ref{fig10} and \ref{fig11}, respectively.
    
From these results, we observe that for $p=1(d=1)$ and $p=2(d=1)$ the optimization stops for $80$ and $43$ iterations, with angles $\theta$ equal to $4.72$ and $9.33$, respectively, while in all of the remaining simulations an angle $\theta<4$ is reached, which is a reasonable value for a numerical solutions \cite{Amstutz2006}. We also notice that, although the cost function converges to a similar minimum value, all the shapes present different solutions with different configurations of features. Regarding the minimization, the lowest cost function values are obtained for $ p=3 (d=1)$ with $4.12$ and $ p=3 (d=3)$ with $4.13$, both with $200$ iterations. While, for $p=2(d=2)$, a slightly higher value of $4.15$ is obtained, but requires only $129$ iterations.
  
In topology optimization, the presence of a local minimum is a well-known challenge and has a dependence on the choice of initial parameters \cite{Sigmund1998}. Small variations in the definition of the initial setting of parameters can lead to different solutions. As a result, even if the initialization is too far from the global minimum, we still can achieve a solution that converges to a local minimum \cite{Allaire2004}. Therefore, while strategies like refining the mesh in a coarse-to-fine approach can partially cure this issue, they do not eliminate it \cite{Allaire2006}. In this example, we notice that the possibility of setting different configurations of polynomial degrees, for the level-set discretization and the solution approximation, does not overcome this phenomenon completely, but depending on the choice of setting, we can obtain a solution that satisfies the condition of having a small $\theta$. Another parameter that can be changed and opens the possibility to find different solutions is the parameter $\gamma$, which limits the size of features in the filtering process, overcoming some noise in the final shape.

%% file: section7_Conclusion.tex
\section{Conclusion}
In this work, we develop an immersed isogeometric approach for a level-set based topology optimization guided only by the topological derivative. The isogeometric approach within this framework provides a seamless geometry update due to a simplified mesh process, defined by a knot vector and control points. It also facilitates straightforward higher-order simulations, yielding results that are comparable to or slightly better than those obtained by the standard approach. In addition, we also treat the elements cut by the level-set function using a quadrature library for implicitly defined geometries to compute a material property used to neglect the contributions from outside the domain $\Omega$, and using a filtering process to smooth the change of material between elements.

This study investigates the impact of using higher-degree basis functions in both the approximation of the solution and the discretization of the level-set function. In this investigation, we observe that being able to perform higher-order simulations can be beneficial, in terms of iterations or the final cost function. However, regarding the level-set discretization, we observe that using linear basis functions yields results comparable to those obtained with higher-degree polynomials. Therefore, although the level-set is continuous, there is a discontinuity in the material property between the material region and the void. Our results indicate that using higher-degree basis functions does not directly imply a better representation of the jump across the interface. We will further investigate this aspect in future research.

%% file: main.bbl
\begin{thebibliography}{34}
\providecommand{\natexlab}[1]{#1}
\providecommand{\url}[1]{\texttt{#1}}
\expandafter\ifx\csname urlstyle\endcsname\relax
  \providecommand{\doi}[1]{doi: #1}\else
  \providecommand{\doi}{doi: \begingroup \urlstyle{rm}\Url}\fi

\bibitem[Allaire and Jouve(2006)]{Allaire2006}
G.~Allaire and F.~Jouve.
\newblock Coupling the level set method and the topological gradient in structural optimization, [in:] \textit{Proceedings {IUTAM} {Symposium} on {Topological} {Design} {Optimization} of {Structures}, {Machines} and {Materials}}, pp. 3-12.
\newblock Dordrecht, 2006.
\newblock \doi{https://doi.org/10.1007/1-4020-4752-5_1}.

\bibitem[Allaire et~al.(2002)Allaire, Jouve, and Toader]{Allaire2002}
G.~Allaire, F.~Jouve, and A.~Toader.
\newblock A level-set method for shape optimization.
\newblock \emph{C. R. Acad. Sci. Paris, Ser. I}, \textbf{334}:\penalty0 1125--1130, 2002.
\newblock \doi{https://doi.org/10.1016/S1631-073X(02)02412-3}.

\bibitem[Allaire et~al.(2004)Allaire, Jouve, and Toader]{Allaire2004}
G.~Allaire, F.~Jouve, and A.~Toader.
\newblock Structural optimization using sensitivity analysis and a level-set method.
\newblock \emph{Journal of Computational Physics}, \textbf{194}:\penalty0 363--393, 2 2004.
\newblock \doi{https://doi.org/10.1016/j.jcp.2003.09.032}.

\bibitem[Aminzadeh and Tavakkoli(2022)]{Aminzadeh2022}
M.~Aminzadeh and S.~M. Tavakkoli.
\newblock A parameter space approach for isogeometrical level set topology optimization.
\newblock \emph{International Journal for Numerical Methods in Engineering}, \textbf{123}:\penalty0 3485--3506, 8 2022.
\newblock \doi{https://doi.org/10.1002/nme.6976}.

\bibitem[Aminzadeh and Tavakkoli(2024)]{Aminzadeh2024}
M.~Aminzadeh and S.~M. Tavakkoli.
\newblock Multiscale topology optimization of structures by using isogeometrical level set approach.
\newblock \emph{Finite Elements in Analysis and Design}, \textbf{235}, 8 2024.
\newblock \doi{https://doi.org/10.1016/j.finel.2024.104167}.

\bibitem[Amstutz(2006)]{Amstutz2006sens}
S.~Amstutz.
\newblock Sensitivity analysis with respect to a local perturbation of the material property.
\newblock \emph{Asymptotic Analysis}, \textbf{49}\penalty0 (1-2):\penalty0 87--108, 2006.
\newblock \doi{https://doi.org/10.3233/ASY-2006-778}.

\bibitem[Amstutz and Andrä(2006)]{Amstutz2006}
S.~Amstutz and H.~Andrä.
\newblock A new algorithm for topology optimization using a level-set method.
\newblock \emph{Journal of Computational Physics}, \textbf{216}:\penalty0 573--588, 8 2006.
\newblock \doi{https://doi.org/10.1016/j.jcp.2005.12.015}.

\bibitem[Bendsøe and Kikuchi(1988)]{Bendsøe1998}
M.~P. Bendsøe and N.~Kikuchi.
\newblock Generating optimal topologies in structural design using a homogenization method.
\newblock \emph{Computer Methods in Applied Mechanics and Engineering}, \textbf{71}:\penalty0 197--224, 1988.
\newblock \doi{https://doi.org/10.1016/0045-7825(88)90086-2}.

\bibitem[Burger et~al.(2004)Burger, Hackl, and Ring]{Burger2004}
M.~Burger, B.~Hackl, and W.~Ring.
\newblock Incorporating topological derivatives into level set methods.
\newblock \emph{Journal of Computational Physics}, \textbf{194}:\penalty0 344--362, 2 2004.
\newblock \doi{https://doi.org/10.1016/j.jcp.2003.09.033}.

\bibitem[de~Prenter et~al.(2023)de~Prenter, Verhoosel, van Brummelen, Larson, and Badia]{Prenter2023}
F.~de~Prenter, C.~V. Verhoosel, E.~H. van Brummelen, M.~G. Larson, and S.~Badia.
\newblock Stability and conditioning of immersed finite element methods: Analysis and remedies.
\newblock \emph{Archives of Computational Methods in Engineering}, \textbf{30}:\penalty0 3617--3656, 7 2023.
\newblock \doi{https://doi.org/10.1007/s11831-023-09913-0}.

\bibitem[Dedè et~al.(2012)Dedè, Borden, and Hughes]{dede2012}
L.~Dedè, M.~J. Borden, and T.~J.R. Hughes.
\newblock Isogeometric analysis for topology optimization with a phase field model.
\newblock \emph{Archives of Computational Methods in Engineering}, \textbf{19}:\penalty0 427--465, 9 2012.
\newblock \doi{https://doi.org/10.1007/s11831-012-9075-z}.

\bibitem[Gangl(2020)]{Gangl2020}
P.~Gangl.
\newblock A multi-material topology optimization algorithm based on the topological derivative.
\newblock \emph{Computer Methods in Applied Mechanics and Engineering}, \textbf{366}, 7 2020.
\newblock \doi{https://doi.org/10.1016/j.cma.2020.113090}.

\bibitem[Gangl and Sturm(2020)]{GanglSturm2020}
P.~Gangl and K.~Sturm.
\newblock A simplified derivation technique of topological derivatives for quasi-linear transmission problems.
\newblock \emph{ESAIM - Control, Optimisation and Calculus of Variations}, \textbf{26}, 2020.
\newblock \doi{https://doi.org/10.1051/cocv/2020035}.

\bibitem[Gangl et~al.(2025)Gangl, Komann, Krenn, and Ulbrich]{Gangl2025}
P.~Gangl, T.~Komann, N.~Krenn, and S.~Ulbrich.
\newblock Robust topology optimization of electric machines using topological derivatives.
\newblock \emph{Arxiv}, 4 2025.
\newblock URL \url{http://arxiv.org/abs/2504.05070}.

\bibitem[Hughes et~al.(2005)Hughes, Cottrell, and Bazilevs]{Hughes2005}
T.J.R. Hughes, J.~A. Cottrell, and Y.~Bazilevs.
\newblock Isogeometric analysis: Cad, finite elements, nurbs, exact geometry and mesh refinement.
\newblock \emph{Computer Methods in Applied Mechanics and Engineering}, \textbf{194}:\penalty0 4135--4195, 10 2005.
\newblock \doi{https://doi.org/10.1016/j.cma.2004.10.008}.

\bibitem[Jahangiry and Tavakkoli(2017)]{Jahangiry2017}
H.~A. Jahangiry and S.~M. Tavakkoli.
\newblock An isogeometrical approach to structural level set topology optimization.
\newblock \emph{Computer Methods in Applied Mechanics and Engineering}, \textbf{319}:\penalty0 240--257, 6 2017.
\newblock \doi{https://doi.org/10.1016/j.cma.2017.02.005}.

\bibitem[Khatibinia et~al.(2020)Khatibinia, Khatibinia, and Roodsarabi]{Khatibinia2020}
M.~Khatibinia, M.~Khatibinia, and M.~Roodsarabi.
\newblock Structural topology optimization based on hybrid of piecewise constant level set method and isogeometric analysis.
\newblock \emph{International Journal of Optimization in Civil Engineering}, \textbf{10}:\penalty0 493--512, 2020.
\newblock URL \url{https://www.researchgate.net/publication/342924423}.

\bibitem[Krenn(2021)]{Krenn2021}
N.~Krenn.
\newblock Multi-material topology optimization subject to pointwise stress constraints for {Additive} {Manufacturing}.
\newblock Master's thesis, Graz University of Technology, 2021.

\bibitem[Ma et~al.(2020)Ma, Zheng, Lei, Zhu, Jin, and Guo]{Ma2020}
B.~Ma, J.~Zheng, G.~Lei, J.~Zhu, P.~Jin, and Y.~Guo.
\newblock Topology optimization of ferromagnetic components in electrical machines.
\newblock \emph{IEEE Transactions on Energy Conversion}, \textbf{35}:\penalty0 786--798, 6 2020.
\newblock \doi{https://doi.org/10.1109/TEC.2019.2960519}.

\bibitem[Osher and Sethian(1988)]{Osher1988}
S.~Osher and J.~A. Sethian.
\newblock Fronts propagating with curvature-dependent speed: Algorithms based on hamilton-jacobi formulations.
\newblock \emph{Journal of Computational Physics}, \textbf{79}:\penalty0 12--49, 1988.
\newblock \doi{https://doi.org/10.1016/0021-9991(88)90002-2}.

\bibitem[Roodsarabi et~al.(2016)Roodsarabi, Khatibinia, and Sarafrazi]{Roodsarabi2016}
M.~Roodsarabi, M.~Khatibinia, and S.~R. Sarafrazi.
\newblock Hybrid of topological derivative-based level set method and isogeometric analysis for structural topology optimization.
\newblock \emph{Steel and Composite Structures}, \textbf{21}:\penalty0 1389--1410, 8 2016.
\newblock \doi{https://doi.org/10.12989/scs.2016.21.6.1389}.

\bibitem[Saye(2015)]{Saye2015}
R.~I. Saye.
\newblock High-order quadrature methods for implicitly defined surfaces and volumes in hyperrectangles.
\newblock \emph{SIAM Journal on Scientific Computing}, \textbf{37}:\penalty0 A993--A1019, 2015.
\newblock \doi{https://doi.org/10.1137/140966290}.

\bibitem[Saye(2022)]{Saye2022}
R.~I. Saye.
\newblock High-order quadrature on multi-component domains implicitly defined by multivariate polynomials.
\newblock \emph{Journal of Computational Physics}, \textbf{448}, 1 2022.
\newblock \doi{10.1016/j.jcp.2021.110720}.

\bibitem[Schillinger and Ruess(2015)]{Schillinger2015}
D.~Schillinger and M.~Ruess.
\newblock The finite cell method: A review in the context of higher-order structural analysis of cad and image-based geometric models.
\newblock \emph{Archives of Computational Methods in Engineering}, \textbf{22}:\penalty0 391--455, 7 2015.
\newblock \doi{https://doi.org/10.1007/s11831-014-9115-y}.

\bibitem[Shojaee et~al.(2012)Shojaee, Mohamadian, and Valizadeh]{Shojaee2012}
S.~Shojaee, M.~Mohamadian, and N.~Valizadeh.
\newblock Composition of isogeometric analysis with level set method for structural topology optimization.
\newblock \emph{International Journal of Optimization in Civil Engineering}, \textbf{2}:\penalty0 47--70, 2012.
\newblock URL \url{https://www.researchgate.net/publication/259593893}.

\bibitem[Sigmund and Maute(2013)]{Sigmund2013}
O.~Sigmund and K.~Maute.
\newblock Topology optimization approaches: A comparative review.
\newblock \emph{Structural and Multidisciplinary Optimization}, \textbf{48}:\penalty0 1031--1055, 12 2013.
\newblock \doi{https://doi.org/10.1007/s00158-013-0978-6}.

\bibitem[Sigmund and Petersson(1998)]{Sigmund1998}
O.~Sigmund and J.~Petersson.
\newblock Numerical instabilities in topology optimization: A survey on procedures dealing with checkerboards, mesh-dependencies and local minima.
\newblock \emph{Structural Optimization}, \textbf{16}:\penalty0 68--75, 1998.
\newblock \doi{https://doi.org/10.1007/BF01214002}.

\bibitem[Teixeira et~al.(2025)Teixeira, Loibl, and Marussig]{Teixeira2025}
G.~H. Teixeira, M.~Loibl, and B.~Marussig.
\newblock Comparison of integration methods for cut elements.
\newblock \emph{ArXiv}, 1 2025.
\newblock \doi{https://doi.org/10.23967/eccomas.2024.098}.
\newblock URL \url{http://arxiv.org/abs/2501.03854}.

\bibitem[Toprak et~al.(2025)Toprak, Loibl, Teixeira, Shiskina, Miao, Kiendl, Marussig, and Kummer]{toprak2025}
T.~Toprak, M.~Loibl, G.~H. Teixeira, I.~Shiskina, C.~Miao, J.~Kiendl, B.~Marussig, and F.~Kummer.
\newblock Employing continuous integration inspired workflows for benchmarking of scientific software -- a use case on numerical cut cell quadrature.
\newblock \emph{ArXiv}, 4 2025.
\newblock URL \url{https://arxiv.org/abs/2503.17192}.

\bibitem[Verzicco(2023)]{Verzicco2023}
R.~Verzicco.
\newblock Immersed boundary methods: Historical perspective and future outlook.
\newblock \emph{Annual Review of Fluid Mechanics}, \textbf{11}:\penalty0 39, 2023.
\newblock \doi{https://doi.org/10.1146/annurev-fluid-120720}.

\bibitem[Vázquez(2016)]{Vázquez2016}
R.~Vázquez.
\newblock A new design for the implementation of isogeometric analysis in octave and matlab: Geopdes 3.0.
\newblock \emph{Computers and Mathematics with Applications}, \textbf{72}:\penalty0 523--554, 8 2016.
\newblock \doi{https://doi.org/10.1016/j.camwa.2016.05.010}.

\bibitem[Wang et~al.(2003)Wang, Wang, and Guo]{Wang2003}
M.~Y. Wang, X.~Wang, and D.~Guo.
\newblock A level set method for structural topology optimization.
\newblock \emph{Computer Methods in Applied Mechanics and Engineering}, \textbf{192}:\penalty0 227--246, 2003.
\newblock \doi{https://doi.org/10.1016/S0045-7825(02)00559-5}.

\bibitem[Wang et~al.(2009)Wang, Zhang, and Liu]{Wang2009}
X.~S. Wang, L.~T. Zhang, and W.~K. Liu.
\newblock On computational issues of immersed finite element methods.
\newblock \emph{Journal of Computational Physics}, \textbf{228}:\penalty0 2535--2551, 4 2009.
\newblock \doi{https://doi.org/10.1016/j.jcp.2008.12.012}.

\bibitem[Wiesheu et~al.(2024)Wiesheu, Komann, Merkel, Schöps, Ulbrich, and Garcia]{Wiesheu2024}
M.~Wiesheu, T.~Komann, M.~Merkel, S.~Schöps, S.~Ulbrich, and I.~Cortes Garcia.
\newblock Combined parameter and shape optimization of electric machines with isogeometric analysis.
\newblock \emph{Optimization and Engineering}, \textbf{26}:\penalty0 1011--1038, 2024.
\newblock \doi{https://doi.org/10.1007/s11081-024-09925-0}.

\end{thebibliography}
